\documentclass[10pt,a4paper]{article}

\usepackage{url}
\usepackage[utf8]{inputenc}
\usepackage[T1]{fontenc}
\usepackage[english]{babel}
\usepackage{a4wide}
\usepackage{float}
\usepackage{graphicx}
\everymath{\displaystyle}
\usepackage{amscd}
\usepackage{amsfonts}
\usepackage{amsmath}
\usepackage{amssymb}
\usepackage{amsthm}
\usepackage{fancyhdr}
\usepackage{latexsym}
\usepackage{mathrsfs}

\newcommand{\cqfd}{\begin{flushright} $\square$ \end{flushright} \end{demo}}
\newcommand{\ensemble}[2]{\left\{ \left. #1 \, \right| \, #2 \right\}}
\newcommand{\Ensemble}[2]{\left\{ #1 \, \left| \, #2 \right. \right\}}
\newcommand{\ec}{\mathbb C}
\newcommand{\ee}{\mathbb E}
\newcommand{\ed}{\mathbb D}
\newcommand{\ea}{\mathbb A}
\newcommand{\FR}{\mathcal{F}}

\DeclareMathOperator{\Ext}{Ext}

\newcommand{\SL}{SL}
\DeclareMathOperator{\Sub}{Sub}
\DeclareMathOperator{\md}{mod}
\DeclareMathOperator{\End}{End}
\DeclareMathOperator{\soc}{soc}
\newcommand{\er}{\mathbb R}
\newcommand{\en}{\mathbb N}
\newcommand{\ez}{\mathbb Z}
\newcommand{\ER}{\mathcal{E}}
\newcommand{\RR}{\mathcal{R}}
\newcommand{\SR}{\mathcal{S}}

\title{Total positivity criteria for partial flag varieties}
\author{Nicolas Chevalier
}
\date{}

\newtheorem{thm}{Theorem}[section]
\newtheorem{prop}[thm]{Proposition}
\newtheorem{lem}[thm]{Lemma}
\newtheorem{problem}[thm]{Problem}

\theoremstyle{definition}
\newtheorem{defi}[thm]{Definition}

\def\gg{\mathfrak{g}}
\def\nn{\mathfrak{n}}
\def\hh{\mathfrak{h}}

\begin{document}

\maketitle

\begin{abstract}
For a simply-connected complex algebraic group $G$ of type
$\ea$, $\ed$, or $\ee$, we prove (see \cite{GLS6}, Conjecture 19.2)
a new family of total positivity criteria for partial flag varieties
$G/P$, where $P$ is a parabolic subgroup of $G$.
\end{abstract}


%
%

\section{Introduction}
\label{sect1}


A matrix $x$ with real entries is called totally positive
if all
its minors
are positive.
These matrices were first studied by I.~Schoenberg~\cite{schoenberg} in the 1930s,
then by F.~Gantmacher and M.~Krein~\cite{GK}, who
showed that the eigenvalues of an $n\times n$ totally positive matrix
are real, positive, and distinct.

G. Lusztig extended this classical subject by introducing
first (in \cite{Lu2}) the totally positive
variety $G_{>0}$
in an arbitrary reductive group $G$, then
(in \cite{Lu1}) the totally positive varieties
$\left( P \setminus G \right) _{>0}$
for any parabolic subgroup $P$ of $G$. Lusztig showed that
$\left( P \setminus G \right) _{>0}$ can be defined by algebraic inequalities
involving the canonical bases.

In 2001, S. Fomin and A. Zelevinsky~\cite{FZ2} introduced
the class of cluster algebras with the purpose of building
a combinatorial framework for studying total positivity in algebraic groups
and canonical bases in quantum groups.
Cluster algebras can be found in several areas of mathematics (for instance combinatorics, Lie theory,
mathematical physics and representation theory of algebras). Other connexions are listed
on Fomin's cluster algebras portal \cite{F}.

C. Geiss, B. Leclerc and J. Schröer have studied cluster algebras associated with Lie groups of type $\ea$, $\ed$, $\ee$,
and have modelled them by categories of modules over the Gelfand-Ponomarev preprojective algebras $\Lambda$
of the same type \cite{GP} (see also \cite{R}).
They have shown \cite{GLS1} that each reachable maximal
rigid $\Lambda$-module can be thought of as a seed of a cluster algebra structure on
$\ec[N]$, the coordinate ring of a maximal unipotent subgroup of $G$
(here, reachable
means that the maximal rigid module is obtained from a distinguished one by a sequence of mutations,
see the end of Section \ref{sect8} for a more precise definition).
They also attached
to each standard parabolic subgroup $P$ of $G$ a certain subcategory $\mathcal{C}_P$ of $\md \Lambda$
and showed that each reachable maximal rigid $\Lambda$-module in $\mathcal{C}_P$ gives a seed for a cluster algebra structure on
$\ec[N_P]$, the coordinate ring of the unipotent radical of $P$.

\begin{problem}[\cite{GLS6}, Conjecture 19.2]
Each basic maximal rigid $\Lambda$-module
in $\mathcal{C}_P$ gives rise to a total positivity criterium for the partial flag variety $P \setminus G$.
\end{problem}

In this note we show (Theorem \ref{resultat1}) that every
reachable basic maximal rigid module in $\mathcal{C}_P$ gives rise to a total positivity criterium.
This leads to a (generally infinite) number of criteria.
These criteria were previously known in the following cases : when $P$ is
a Borel subgroup, that is for the total flag variety (Berenstein,
Fomin and Zelevinsky~\cite{BFZ}, and \cite{FZ1}), and
when $P \setminus G$ is a type $\ea$ grassmannian (Scott
\cite{scott}). In all other cases, for example
for partial flag varieties in type $\ea$, or
for Grassmannians in type $\ed$ and $\ee$, these criteria are new.

In fact, the proof of Theorem \ref{resultat1} turns out to be rather easy if one suitably combines
information coming from several sources \cite{FZ1,GLS1,GLS3,GLS4}.
The main idea is to use the algorithm of \cite{GLS4}, \S 13.1, to relate
maximal rigid modules of \cite{GLS3} with positivity criteria of \cite{FZ1}
(see below, Section \ref{sect7}).

Note that the criteria given here are of the form $\varphi_M(x) > 0$ for some regular functions $\varphi_M$
on $N$ attached to certain rigid $\Lambda$-modules $M$ (their definition will be recalled
in Section \ref{sect3}). Geiss, Leclerc and Schröer showed that these functions belong to the dual
of Lusztig's semicanonical basis \cite{Lu3} of $U(\mathfrak{n})$ (where $\mathfrak{n} = \mathrm{Lie}(N)$).
Since Lusztig expressed total positivity in terms of the canonical basis of $U(\mathfrak{n})$,
this gives some supporting evidence for the conjecture of \cite{GLS1}, stating that the functions $\varphi_M$
for rigid $M$ belong at the same time to Lusztig's dual canonical and dual semicanonical basis.


\section{Flag varieties and their totally positive part}
\label{sect2}


Let $\gg$ be a simple complex Lie algebra of rank~$n$ of type $\ea$, $\ed$, or $\ee$,
with the Cartan decomposition $\gg = \nn \oplus \hh \oplus \nn^-$.
Let $e_i \in \nn, h_i \in \hh, f_i \in \nn^- \,$, for $i \in I=\left\lbrace 1, \ldots, n \right\rbrace $ be the
Chevaley generators of~$\gg$, and let $A = (a_{ij})$ be the Cartan matrix.
Thus $a_{ij} = \alpha_j (h_i)$, where 
$\alpha_1, \ldots, \alpha_n \in \hh^*$ are the simple roots of~$\gg$.
Let $R$ denote the root system of $g$.
Let $G$ be a simply connected complex Lie group with
the Lie algebra $\gg$.
Let $N^-$, $H$ and $N$ be closed subgroups of $G$ with Lie algebras
$\nn^-$, $\hh$ and $\nn$, respectively.
Thus $H$ is a maximal torus,
and $N$ and $N^-$ are two opposite maximal unipotent subgroups of $G$.
Let $B^- = HN^-$ and $B = HN$ be the corresponding pair of opposite
Borel subgroups (thus we have $B \cap B^- = H$).

For $i \in I$ and $t \in \ec$, we write
$$ x_i (t) = \exp(t e_i)\ , $$
$$y_{i} (t) = \exp(t f_i) \ ,$$ so that $t \mapsto x_i (t)$ (resp.\ $t \mapsto y_{i} (t)$)
is a one-parameter subgroup in $N$ (resp.\ in $N^-$) and we denote it by $U_{\alpha_i}$
(resp. by $U_{-\alpha_i}$). The one-parameter root subgroup
associated to $\alpha \in R$ is also denoted by $U_{\alpha}$.

The Weyl group $W$ of $G$ is defined by $W = {\rm Norm}_G (H)/H$.
The group $W$ is a Coxeter group with Coxeter generators the simple reflections
$s_1, \ldots, s_n$.
A reduced word for $w \in W$ is a sequence of indices 
$\textbf{i} = (i_1, \ldots, i_m)$ of shortest possible length
such that $w = w(\textbf{i}) = s_{i_1} \cdots s_{i_m}\,$.
The number~$m$ is denoted by $\ell(w)$ and
is called the length of~$w$.
The set of reduced words for~$w$ will be denoted by~$\mathcal{R}(w)$.
The Weyl group $W$ has the unique element $w_0$ of maximal length equal to
$r = \ell (w_0)$, which is also the dimension of the affine space $N$.

For a fixed subset $K \subset I$, we let $B_K$ (resp. $B_K^-$) be the standard parabolic subgroup of $G$ generated
by $B$ and the $\{y_k(t)\}_{k \in K}$ (resp. by $B^-$ and the $\{x_k(t)\}_{k \in K}$).
We let $N_K$ (resp. $N_K^{-}$) be the unipotent radical of $B_K$ (resp. $B_K^-$).
Let $X_K := B_K^{-} \setminus G$ and $\pi_K := G \twoheadrightarrow X_K$ be the canonical
projection. The set $X_K$ is a projective variety called a \emph{partial flag variety}.

\begin{lem}[\cite{Lu2}, \S 2.7 and 2.10]
Let $w \in W$ and let $\textbf{i} \in \mathcal{R}(w)$. Then the image of :
$$ \begin{array}{ccc}
(\er_{>0})^k & \to & N \\
(t_1,t_2,\hdots,t_k) & \mapsto & x_{i_1}(t_1)x_{i_2}(t_2)\hdots x_{i_k}(t_k)
   \end{array} $$
does not depend on the choice of $\textbf{i} \in \mathcal{R}(w)$. We denote it by
$N_{>0}^w$. When $w=w_0$, we write $N_{>0}:=N_{>0}^{w_0}$.
\end{lem}

We are now able to define $(X_K)_{>0}$.

\begin{defi}[\cite{Lu1}, \S 1.5]
The totally positive part of $X_K$ is :
$$ (X_K)_{>0} = \pi_K(N_{>0}) . $$
\end{defi}


\section{The preprojective algebra $\Lambda$} 
\label{sect3}


Let $\overline{Q}$ denote the quiver obtained
from the Dynkin diagram of $\gg$ by replacing
every edge by a pair $(a,a^*)$
of opposite arrows.
Let
\[
 \rho = \sum (aa^* - a^*a)
\]
be the element of the path algebra $\ec\overline{Q}$ of $\overline{Q}$,
where the sum is over all pairs of opposite arrows.
Following \cite{GP,R}, we define the \emph{preprojective 
algebra} $\Lambda$ as the quotient of $\ec\overline{Q}$ by 
the two-sided ideal generated by $\rho$.
This is a finite-dimensional selfinjective algebra,
with infinitely many isomorphism classes of indecomposable
modules, except if $\gg$ has type $\ea_n$ with $n\leqslant 4$.
The category $\md \Lambda$ has the following
important symmetry property (see \cite{GLS2}) :
$$ \forall M,N \in \md \Lambda, ~~~~~~~~~~~~ \dim \Ext_{\Lambda}^1
(M,N) = \dim \Ext_{\Lambda}^1 (N,M). $$
A $\Lambda$-module $M$ is said to be \textit{rigid} when $\Ext_{\Lambda}^1(M,M)=0$.
Let $S_i \ (1\le i\le n)$ be the one-dimensional $\Lambda$-modules
attached to the vertices $i$ of $\overline{Q}$. We let $P_i$ and $Q_i$ be respectively
the projective cover and the injective hull of $S_i$.

There is a close relationship between
$\Lambda$ and $\ec[N]$.
Geiss, Leclerc and Schr\"oer \cite{GLS1,GLS2} have attached to every object $M$ in $\md\Lambda$
a polynomial function $\varphi_M$ on $N$.
These functions may be defined as follows (see \cite{GLS1}, Lemma 9.1).
For a multi-integer
$\textbf{a} = (a_1,a_2, \hdots, a_k) \in \ez^k_{\geqslant 0}$,
for $\textbf{t}=(t_1,\hdots,t_k) \in \ec^k$
and for a multi-index $\textbf{i}=(i_1,i_2, \hdots , i_k) \in I^k$, we write :
$$ \begin{array}{rcl}
\textbf{a} ! &:=& a_1 ! a_2 ! \hdots a_k !,\\
\textbf{t}^{\textbf{a}} &:=& t_1^{a_1} t_2^{a_2} \hdots t_k^{a_k},\\
x_{\textbf{i}}(\textbf{t}) &:=& x_{i_1}(t_1) \hdots x_{i_k}(t_k) \in N,\\
\textbf{i}^{\textbf{a}} &:=&
(\underbrace{i_1,i_1, \hdots ,i_1}_{a_1 \mbox{\scriptsize{ times}}}, \hdots ,
\underbrace{i_k,i_k, \hdots ,i_k}_{a_k \mbox{\scriptsize{ times}}}). \end{array}$$
For a $\Lambda$-module $M$, we denote by :
$$ \mathfrak{f}:= \big( \{0\}=M_0 \subset M_1 \subset \hdots \subset M_d=M \big)$$
a composition series of $M$, that is a flag of sub-$\Lambda$-modules of $M$ where all
consecutive quotients are simple : there exists $i_k \in I$ such that
$M_k / M_{k-1}
\cong S_{i_k}$. We call $\textbf{i}:=(i_1,i_2, \hdots , i_d) \in I^d$ the \emph{type} of $ \mathfrak{f}$.
We denote by $\Phi_{\textbf{i},M}$ the projective variety of flags of $M$ whose type is
$\textbf{i}$, and by $\chi_{\textbf{i},M}:=\chi(\Phi_{\textbf{i},M})$ its Euler characteristic.
With this notation, we can state the following lemma :
\begin{lem}[\cite{GLS1}, Lemma 9.1]\label{formule}
Let $\textbf{i}=(i_1,i_2, \hdots, i_k) \in I^k$ and $M \in \md \Lambda$. Then :
$$ \varphi_M \left( x_{\textbf{i}}(\textbf{t}) \right)= \sum_{\textbf{a} \in \en^k}
\chi_{\textbf{i}^{\textbf{a}},M} \dfrac{\textbf{t}^{\textbf{a}}}{\textbf{a}!}.$$
\end{lem}


\section{Total positivity criteria for $N$}
\label{sect4}


In order to know whether an element $n \in N$ lies in $N_{>0}$, Fomin and Zelevinsky
\cite{FZ1} gave a positivity criterium for each $\textbf{i} \in \RR(w_0)$. First
we construct $\underline{w}$ and $\underline{\underline{w}}$, two representatives
of $w \in W$ in $G$. Let $\phi_i : \SL_2(\mathbb{C}) \to  G$ be the group
homomorphism defined by :
$$ \phi_i \left( \begin{array}{cc} 1 & t \\ 0 & 1 \end{array} \right) = x_i(t),$$
$$ \phi_i \left( \begin{array}{cc} 1 & 0 \\ t & 1 \end{array} \right) = y_i(t).$$
Define :
$$\underline{s_i} = \phi_i \left( \begin{array}{cc} 0 & -1 \\ 1 & 0 \end{array} \right),$$
$$\underline{\underline{s_i}} = \phi_i \left( \begin{array}{cc} 0 & 1 \\ -1 & 0 \end{array} \right).$$
where $s_i$ is a simple reflection in $W$. If $\textbf{i} \in \RR(w_0)$, then the following elements
are well defined in $G$ :
$$ \begin{array}{ccc}
\underline{w} & = & \underline{s_{i_1}} ~ \underline{s_{i_2}} \hdots \underline{s_{i_l}}, \\
\underline{\underline{w}} & = & \underline{\underline{s_{i_1}}} ~ \underline{\underline{s_{i_2}}}
\hdots \underline{\underline{s_{i_l}}}.
   \end{array}$$

Recall that the weight lattice is the set of all weights $\gamma \in \hh^*$
such that $\gamma (h_i) \in \ez$ for all $i$.
It has a $\ez$-basis formed by the fundamental weights
$\varpi_1, \ldots, \varpi_n$ defined by $\varpi_i (h_j) = \delta_{ij}$.
Every such weight $\gamma$ gives rise to a multiplicative character
$a \mapsto a^\gamma$ of the maximal torus $H$; this character
is given by $\exp (h)^\gamma = e^{\gamma (h)}  \,\, (h \in \hh)$.

Let $G_0$ be the subset of $G$ whose elements admit a gaussian reduction, that is for all
$g \in G_0$, one can write $g=[g]_-[g]_0[g]_+$ with $[g]_- \in N^{-}$, $[g]_0 \in H$, and $[g]_+ \in N$.
Let $(\Delta^{\varpi_i})_{i \in I}$ be the regular functions on $G$ which satisfy the following condition :
for all $g \in G_0$,
$\Delta^{\varpi_i}(g)=\Delta^{\varpi_i}([g]_0)$, and if $g \in H$, then
$\Delta^{\varpi_i}(g)=g^{\varpi_i}$.

\begin{defi}[\cite{FZ1}, Definition 1.4]
Let $u,v \in W$ and $i \in I$. For $x \in G$, put :
$$ \Delta_{u (\varpi_i) , v (\varpi_i)}(x):=\Delta^{\varpi_i}(\underline{\underline{u^{-1}}}x \underline{v}).$$
These functions are called \emph{generalized minors} of $x$. We denote by
$D_{u (\varpi_i) , v (\varpi_i)}$ the restriction of $\Delta_{u (\varpi_i) , v (\varpi_i)}$ to $N$.
\end{defi}

In type $\ea$, $\Delta_{u (\varpi_i) , v (\varpi_i)}(x)$ is nothing else but
the classical minor of $x$ of size $i$ corresponding to the submatrix with row set $\{u(1),u(2), \hdots, u(i) \}$
and column set $\{ v(1), v(2), \hdots, v(i) \}$ (where we identify $u$ and $v$ to permutations of $\{ 1,\hdots,n \}$).

\begin{defi}\label{t_l}     
Let $\textbf{i}=(i_1, i_2, \hdots, i_r) \in \mathcal{R}(w_0)$. For all
$l \in I$, put :
$$ t_l := \max \ensemble{t \leqslant r}{i_t=l}.$$
\end{defi}

This is the right-most index of $\textbf{i}$ equal to $l$.
The next lemma is straightforward.

\begin{lem}
Let $\textbf{i}=(i_1, i_2, \hdots, i_r) \in \RR(w_0)$ and $l \in I$.
We have :
$$ \begin{array}{ccrcl}
\forall x \in G, & ~~~~~~~~~~~~ & \Delta_{\varpi_{i_{t_l}},s_{i_r}s_{i_{r-1}} \hdots
s_{i_{{t_l}+1}} (\varpi_{i_{t_l}}) } (x) & = & \Delta_{\varpi_l, \varpi_l} (x) , \\
\forall n \in N, && D_{\varpi_{i_{t_l}},s_{i_r}s_{i_{r-1}} \hdots
s_{i_{{t_l}+1}} (\varpi_{i_{t_l}}) } (n) & = & 1.
   \end{array}$$
\end{lem}

The minors $\Delta_{\varpi_{i_{t_l}},s_{i_r}s_{i_{r-1}} \hdots
s_{i_{{t_l}+1}} (\varpi_{i_{t_l}}) }$ in the preceding lemma are precisely those which
appear in \cite{FZ1}, formulas (1.18) and (1.23), where $u=e$, $v=w_0$,
$\textbf{i}=(i_r, \hdots, i_1)$ is written backwards, and $k=t_l$ for some $l \in I$.

\begin{defi}\label{cluster}
For $\textbf{i}=(i_1, i_2, \hdots, i_r) \in \RR(w_0)$, define
$$e(\textbf{i}):=\ensemble{m}{1 \leqslant m \leqslant r \mbox{ and } \forall l \in I, \, m \neq t_l}$$
to be the subset of $\{ 1, \hdots, r \}$ where we remove all the $t_l$, for $l \in I$, and let :
$$ \FR(\textbf{i}) := \ensemble{\Delta_{\varpi_i,w_0(\varpi_i)}}{1 \leqslant i \leqslant n}
\cup \ensemble{\Delta_{\varpi_{i_k}, s_{i_r}s_{i_{r-1}} \hdots s_{i_{k+1}}
(\varpi_{i_k})}}{k \in e(\textbf{i})}. $$
\end{defi}

In \cite{FZ1} formula (1.23), the set $F(\textbf{i})$ is equal to the union of $ \FR(\textbf{i}) $
and $\Delta_{\varpi_{i_{t_l}},s_{i_r}s_{i_{r-1}} \hdots
s_{i_{{t_l}+1}} (\varpi_{i_{t_l}}) }$ (for $l \in I$). But
the restriction to $N$ of those last generalized minors
being equal to $1$ (hence positive),
we will not need them.

Define $G^{u,v} := B \underline{u} B \cap B^{-} \underline{v} B^{-}$, the intersection
of the two opposite Schubert cells.
In \cite{Lu2}, Lusztig has defined its positive part $G^{u,v}_{>0}$,
and in \cite{FZ1} Theorem 1.11, Fomin and Zelevinsky have given parametrizations of $G^{u,v}_{>0}$.
Taking into account that $N_{>0}=N \cap G^{e,w_0}_{>0}$, we can state :

\begin{thm}[\cite{FZ1}, Theorem 1.11] \label{Chamber_Ansatz}
Let $\textbf{i}=(i_1, i_2, \hdots, i_r) \in \RR(w_0)$. The map
$N \to \ec^r$ given by $n \mapsto \left( \Delta (n) ,\, \Delta \in \FR(\textbf{i}) \right)$
restricts to a bijection $N_{>0} \stackrel{\sim}{\rightarrow} \er^r_{>0}$.
\end{thm}


\section{Operations on $\Lambda$-modules}
\label{sect5}


\begin{defi}[\cite{GLS3}, \S 5.1]\label{soc1}
Let $M \in \md \Lambda$ and $i \in I$. Let $m_i^{\dagger}(M)$ be the
multiplicity of the simple $S_i$ in the socle of $M$. We put :
$$\ER_i^{\dagger}(M) := M / S_i^{\oplus m_i^{\dagger}(M)}.$$
\end{defi}

We now recall from \cite{GLS4} the definition of $\soc_{(i_k,\hdots,i_1)}(X)$. Note that
$S_i^{\oplus m_i^{\dagger}(M)} = \soc_{(i)}(M)$.

\begin{defi}[\cite{GLS4} \S 2.4]\label{soc2}
For a $\Lambda$-module $X$ and an index $j$, $1 \leqslant j \leqslant n$, we define
$\soc_{(j)}(X):=\soc_{S_j}(X)$ to be the sum of all submodules $U$ of $X$ which are isomorphic
to $S_j$. For $(j_1,\hdots,j_t) \in I^t$,
there is a unique chain
$$ 0 = X_0 \subset X_1 \subset \hdots \subset X_t \subset X$$
of submodules of $X$ such that $X_p / X_{p-1} = \soc_{(j_p)} \left( X / X_{p-1} \right)$.
Define $\soc_{(j_1, \hdots, j_t)}(X):=X_t$.
\end{defi}

The following lemma is clear.

\begin{lem}\label{lemounet}
For a $\Lambda$-module $X$ and for $(j_1,\hdots,j_t) \in I^t$,
we have :
$$\ER_{j_t}^{\dagger} \hdots \ER_{j_1}^{\dagger}(X) = X / \soc_{(j_1, \hdots, j_t)}(X).$$
\end{lem}

The additive functor $\ER_i^{\dagger}$ satisfies
some relations related to braid relations (see \cite{GLS3}, Proposition 5.1). In particular,
if $\textbf{i}=(i_1,i_2, \hdots i_{\ell(w)}) \in \RR(w)$ for $w \in W$, then
$\ER_{i_1}^{\dagger} \ER_{i_2}^{\dagger} \hdots \ER_{i_{\ell(w)}}^{\dagger}$ does not depend on
$\textbf{i}$.

Let $\SR$ be the self-duality of $\md \Lambda$ introduced in
\cite{GLS5}, \S 1.7.
The formula : $$ \ER_i^{\dagger} = \SR \ER_i \SR,$$ defines another additive functor which satisfies the same properties.
These functors allow us to express generalized minors as $\varphi-$functions :

\begin{prop}[\cite{GLS5}, Lemma 5.4]\label{form_phi}
Let $u,v \in W$. Then :
$$ \varphi_{\ER_u^{\dagger} \ER_v(Q_i)} = D_{u(\varpi_i),vw_0(\varpi_i)}.$$
\end{prop}


\section{Maximal rigid modules and their mutations}
\label{Gamma} \label{sect7b}


Recall that $r$ is the dimension of the affine space $N$.
This is also the number of elements of every cluster of $\ec[N]$
(if we include the frozen variables).
Geiss and Schr\"oer have shown \cite{GS} that 
the number of pairwise non-isomorphic indecomposable direct summands
of a rigid $\Lambda$-module is bounded above by $r$.
A rigid module with $r$ non-isomorphic indecomposable summands 
is called \emph{maximal}.

Let $T=T_1\oplus \cdots \oplus T_r$ be a maximal rigid module in $\md \Lambda$,
where every $T_i$ is indecomposable. Define $B = \End_\Lambda T$, a basic
finite-dimensional algebra with simple modules $b_i\ (1\le i\le r)$.
Denote by $\Gamma(T)$ the Gabriel-quiver of $B$, that is, the quiver with
vertex set $\{1,\ldots,r\}$ and $d_{ij}$ arrows
from $i$ to $j$, where $d_{ij}=\dim \Ext^1_B(b_i,b_j)$.

Define $\Sigma(T) = ((\varphi_{T_1},\ldots,\varphi_{T_r}),\,\Gamma(T))$. This
$\Sigma(T)$ will play the role of an initial seed for a geometric cluster algebra
structure on $\ec[N]$. Here, by geometric, we mean that the quiver $\Gamma(T)$ encodes all the information for the
mutation process, or equivalently, that the Fomin-Zelevinsky's mutation matrix is skew-symmetric with associated
quiver $\Gamma(T)$.

\begin{thm}[\cite{GLS1}]
Let $T=T_1\oplus \cdots \oplus T_r$ be a maximal rigid $\Lambda$-module.
Let $T_k$ be a non-projective indecomposable summand of $T$.
There exists a unique indecomposable module $T_k^*\not\cong T_k$
such that $(T/T_k)\oplus T_k^*$ is maximal rigid.
\end{thm}

The maximal rigid module $(T/T_k)\oplus T_k^*$
is called the \emph{mutation of $T$ in direction $k$},
and is denoted by $\mu_k(T)$.

\begin{thm}[\cite{GLS1}]\label{thGLS1}
We have $\Sigma(\mu_k(T)) = \mu_k(\Sigma(T))$, where in the right-hand side
$\mu_k$ stands for the Fomin-Zelevinsky seed mutation for the
cluster algebra structure on $\ec[N]$.
\end{thm}


\section{Construction of some maximal rigid $\Lambda$-modules}
\label{sect6} \label{sect7}


In this section, we will recall the definition of the maximal rigid $\Lambda$-modules
$T_{\textbf{i}}^{\dagger}$ of \cite{GLS3} \S 5.3 and $V_{\textbf{i}}$
of \cite{GLS4} \S 9.8, and we will see the relations
between them and the set $\FR(\textbf{i})$ (see Definition \ref{cluster}).

\subsection{} \label{U_i}
Let $\textbf{i}=(i_1,i_2, \hdots, i_r) \in \RR(w_0)$.
If $m \in -I=\{-n,\hdots,-1\}$, then we let
$M_m:=Q_{-m}$ be the indecomposable injective rigid module. If $m \in e(\textbf{i})$, then
we set $M_m:=\ER_{i_1}^{\dagger} \hdots \ER_{i_m}^{\dagger}(Q_{i_m})$. Define :
$$T_{\textbf{i}}^{\dagger}:=\bigoplus_{m \in -I \cup e(\textbf{i})} M_m.$$
Note that this $\Lambda$-module coincide with the maximal rigid $\Lambda$-module $T_{\textbf{i}}^{\dagger}$
defined in \cite{GLS3} \S 5.3 (see also the proof of \cite{GLS3}, proposition 6.1).

\subsection{} \label{V_i}
Let $\textbf{i}=(i_r,\hdots,i_1) \in \RR(w_0)$ (beware that we reverse the order of the indices here).
Following \cite{GLS4} \S 9.8, we put
$V_k:=\soc_{(i_k,\hdots,i_1)}(Q_{i_k})$, $V_{\textbf{i}} = \bigoplus_{k=1}^r V_k$
(see Definition \ref{soc2} and \cite{GLS4} \S 2.4 ; here we use our notation $Q_j$ for the injective modules instead of
the notation $\hat{I_j}$ in \cite{GLS4}).
Let $T_{\textbf{i}}:= \bigoplus_{k=1}^r \left( Q_{i_k} / V_{k^-} \right)$, where
$k^{-}:=\max \ensemble{-k,1 \leqslant s \leqslant k-1}{i_s=i_k}$
(here we change a little bit the definition of $k^{-}$ for the convenience of the proofs,
but this change has no real impact on the definitions of the module $T_{\textbf{i}}$).
Both $V_{\textbf{i}}$ and $T_{\textbf{i}}$ are $\Lambda$-modules.

\begin{thm}[\cite{GLS1}]\label{thGLS2}
Let $\textbf{i} \in \RR(w_0)$. The map $T \mapsto \Sigma(T)$ gives a one-to-one correspondence
between the maximal rigid $\Lambda$-modules in the mutation class of $V_{\textbf{i}}$ and the
clusters of $\ec[N]$.
\end{thm}

This theorem (together with Theorem \ref{thGLS1})
allows to lift to $\md \Lambda$ the geometric cluster algebra structure on $\ec[N]$.

\subsection{} \label{lienU_iV_i}

For a multi-index ${\textbf{i}}=(i_1,\ldots,i_l)$, we let
$m(\textbf{i}):=(i_l,\hdots,i_2,i_1)$ be the mirror image of $\textbf{i}$.

\begin{lem}\label{lien}
Let $\textbf{i}=(i_1,\hdots,i_r) \in \RR(w_0)$.
The modules $T_{\textbf{i}}^{\dagger}$ of \cite{GLS3} \S 5.3 and
$ T_{m(\textbf{i})}$ of \cite{GLS4} \S 9.8. coincide.
\end{lem}

\begin{demo}
We have that $ T_{m(\textbf{i})} = \bigoplus_{k=1}^r Q_{i_k} / \soc_{(i_{k^-},\hdots,i_1)}(Q_{i_k})$.
But, by Lemma \ref{lemounet},
$$Q_{i_k} / \soc_{(i_{k^-},\hdots,i_1)}(Q_{i_k}) = \ER_{i_1}^{\dagger} \hdots
\ER_{i_{k^-}}^{\dagger}(Q_{i_k}) = M_{i_{k^-}},$$ hence
$T_{m(\textbf{i})} = \bigoplus_{k=1}^r M_{i_{k^-}}$.
It only remains to show that the sets $-I \cup e(\textbf{i})$ and $ \ensemble{k^-}{1 \leqslant k \leqslant r}$ coincide.
Indeed, they both have cardinality $r$,
so it is enough to show that $-I \cup e(\textbf{i})$
is contained in $\ensemble{k^-}{1 \leqslant k \leqslant r}$.
Now, if $j \in I$ and $s$ is the smallest index $k$ in $\textbf{i}$ such that $i_k=j$ (such an $s$ exists because
$\textbf{i} \in \RR(w_0)$), then $s^{-}=-j$, hence $-I \subset \ensemble{k^-}{1 \leqslant k \leqslant r}$. Finally,
if $j \in e(\textbf{i})$, it means that $j \neq t_l$ for all $l \in I$, hence there exists some index
$s$ (which is the smallest index $k>j$ such that $i_k=i_j$) such that $s^{-}=j$, hence $e(\textbf{i})
\subset \ensemble{k^-}{1 \leqslant k \leqslant r}$ and we are done.
\cqfd

For an illustration of Lemma \ref{lien} in type $\ed_4$, see Example \ref{sect12}.

\begin{prop} \label{algorithme}
Let $\textbf{i} \in \RR(w_0)$. Then
there exists a sequence of mutations in $\md \Lambda$ which begins with $T_{\textbf{i}}^{\dagger}$
and ends at $V_{m(\textbf{i})}$.
\end{prop}

\begin{demo}
In \cite{GLS4}, \S 13.1,  an explicit sequence of mutations between $V_{m(\textbf{i})}$
and $T_{m(\textbf{i})}$ is described algorithmically. The lemma then follows from Lemma \ref{lien}.
\cqfd

\subsection{}

We end this section by relating the module $V_{\textbf{i}}$ of \S \ref{V_i} to the positivity
criterium of Theorem \ref{Chamber_Ansatz}.

\begin{lem}\label{lienV}
Let $\textbf{i}=(i_1,\hdots,i_r) \in \RR(w_0)$, and $V_{\textbf{i}}=V_1 \oplus \cdots \oplus V_r$
as above. Then the sets $\left\{ \varphi_{V_1}, \hdots, \varphi_{V_r} \right\}$ and
$\tilde{\FR}(\textbf{i})$ coincide, where $\tilde{\FR}(\textbf{i})$ is the set
of the restrictions to $N$ of the generalized minors of $\FR(\textbf{i})$.
\end{lem}

\begin{demo}
By \cite{GLS4}, proposition 9.1, we have for $1 \leqslant k \leqslant r$ :
$$ \varphi_{V_{r-k+1}} = D_{\varpi_{i_k},s_{i_r}\hdots s_{i_k}(\varpi_{i_k})}.$$
Note that, in contrast with \cite{GLS4} \S 9, our $\textbf{i}$ is
written backwards here. The two sets $\tilde{\FR}(\textbf{i})$ and
$\left\{ \varphi_{V_1}, \hdots, \varphi_{V_r} \right\}$ have same
cardinality $r$, thus we will only prove that the set $\tilde{\FR}(\textbf{i})$
is included in $\left\{ \varphi_{V_1}, \hdots, \varphi_{V_r} \right\}$.
First, let $j \in I$, then $D_{\varpi_j,w_0(\varpi_j)}=D_{\varpi_j,s_{i_r}\hdots s_{i_1}(\varpi_j)}$,
and if $t$ stands for the smallest index $k$ such that $i_k=j$, we have
$s_{i_{t-1}} \hdots s_{i_1}(\varpi_j)=\varpi_j$, hence
$D_{\varpi_j,w_0(\varpi_j)}=D_{ \varpi_{i_t} , s_{i_r} \hdots s_{i_t} (\varpi_{i_t}) } $,
which prove that $\ensemble{D_{\varpi_j,w_0(\varpi_j)}}{j \in I}$ is included
in $\left\{ \varphi_{V_1}, \hdots, \varphi_{V_r} \right\}$.

Next, let $k \in e(\textbf{i})$. Hence there exists an index $k^{+}$ which is
the smallest index $j>k$ such that $i_j=i_k$. Then $D_{\varpi_{i_k}, s_{i_r}s_{i_{r-1}} \hdots s_{i_{k+1}}
(\varpi_{i_k})}=D_{\varpi_{i_{k^{+}}}, s_{i_r}s_{i_{r-1}} \hdots s_{i_{k^{+}}}(\varpi_{i_{k^{+}}})}$,
hence :
$$\ensemble{D_{\varpi_{i_k}, s_{i_r}s_{i_{r-1}} \hdots s_{i_{k+1}}
(\varpi_{i_k})}}{k \in e(\textbf{i})} \subset \left\{ \varphi_{V_1}, \hdots, \varphi_{V_r} \right\},$$
thus $\tilde{\FR}(\textbf{i}) \subset \left\{ \varphi_{V_1}, \hdots, \varphi_{V_r} \right\}$.
\cqfd

\section{The maximal rigid object $U_{\textbf{i}}^{(K)}$ of $\Sub Q_J$} \label{T_i}

Let $K$ be a fixed subset of $I$ as in Section \ref{sect2}. We denote by $W_K$ the subgroup of $W$ generated by the
$\{s_k\}_{k \in K}$, and $w_0^K$ the element of $W_K$ of maximal length equal to
$r_K = \ell (w_0^K)$. We let
$\RR (w_0,K)$ be the set of reduced words for $w_0$ that are adapted to $K$, that is
if $\textbf{i} = (i_1, \ldots, i_r)$ is in $\RR (w_0,K)$, then
$w(\textbf{i}) = w_0$ and $w(i_1, \ldots, i_{r_K}) = w_0^K$.

For $J := I-K$, we write $Q_J := \bigoplus_{j \in J} Q_j$ and $\Sub Q_J$
is the full subcategory of $\md \Lambda$ whose objects are isomorphic to a submodule of a sum of a finite
number of copies of $Q_J$. Let $\textbf{i} \in \RR(w_0,K)$.
Following \cite{GLS3}, \S 9.3, we construct the object $U_{\textbf{i}}^{(K)}$ of $\Sub Q_J$.
For $k \in K$, let $t_k^{(K)}:=\max \ensemble{t \leqslant r_K}{i_t=k}$, and
for $j \in J$, let $t_j^{(K)}:=-j$. Now define :
$$ \begin{array}{rcl}
I_K & := & \ensemble{t_i^{(K)}}{i \in I} \\~\\
e_K(\textbf{i}) & := & \Ensemble{m}{r_K < m \leqslant r \mbox{ and } m \in e(\textbf{i})} \\~\\
U_{\textbf{i}}^{(K)} & := & \bigoplus_{m \in I_K \cup e_K(\textbf{i})} M_m.
\end{array}$$
It is proved in \cite{GLS3}, Proposition 7.3, \S 9.2 and 9.3, that $U_{\textbf{i}}^{(K)}$
is maximal rigid in $\Sub Q_J$.

\begin{lem}\label{invariance}
Let $X \in \Sub Q_J$. Then $\varphi_X$ is $(B_K^- \cap N)$-invariant. This
is true in particular when $X$ is a direct summand of $U_{\textbf{i}}^{(K)}$.
\end{lem}

\begin{demo}
Fix a $\Lambda$-module $X \in \Sub Q_J$.
The group $B_K^- \cap N$ is the subgroup of $N$ generated
by the $x_i(t)$ for $i \in K$ and $t \in \ec$. Hence we have to prove
that for each $n \in N$, each $k \in K$ and each $t \in \ec$, $\varphi_X(x_k(t)n)
= \varphi_X(n)$. Write $n = x_{i_1}(t_1) x_{i_2}(t_2) \hdots x_{i_l}(t_l)$. Then,
according to Lemma \ref{formule} :
$$ \varphi_X(x_k(t)n) = \varphi_X(x_k(t) x_{i_1}(t_1) x_{i_2}(t_2) \hdots x_{i_l}(t_l))
= \sum_{\textbf{a} \in \en^{l+1}}
\chi_{\textbf{i}^{\textbf{a}},X} \dfrac{\textbf{t}^{\textbf{a}}}{\textbf{a}!},$$
where $\textbf{i} = (k,i_1,i_2,\hdots,i_l)$ and $\textbf{t}^{\textbf{a}}
= t^{a_0} t_1^{a_1} t_2^{a_2} \hdots t_l^{a_l}$. Recall (see Section \ref{sect3}) that $\chi_{\textbf{i}^{\textbf{a}},X}$
is the Euler characteristic of the variety of all composition series of $X$ :
$$ \mathfrak{f}:= \big( \{0\}=X_{-1} \subset X_0 \subset \hdots \subset X_l=X \big)$$
where the $a_0$ first quotients $X_l / X_{l-1}$ are isomorphic to $S_k$ (and $k \in K$),
the $a_1$ next quotients are isomorphic to $S_{i_1}$, etc.
But $X$ belongs to $\Sub Q_j$ and hence has no simple $S_k$ in its socle. It forces
$a_0$ to be $0$ in the preceding formula, and we get :
$$ \varphi_X(x_k(t)n)
= \sum_{\textbf{b} \in \en^{l}} \chi_{\textbf{j}^{\textbf{b}},X} \dfrac{t_1^{b_1} \hdots t_l^{b_l}}{b_1!\hdots b_l!},$$
where
$\textbf{j}=(i_1,\hdots,i_l)$. But :
$$ \varphi_X(n) = \sum_{\textbf{b} \in \en^{l}} \chi_{\textbf{j}^{\textbf{b}},X} \dfrac{t_1^{b_1} \hdots t_l^{b_l}}{b_1!\hdots b_l!}$$
thanks to Lemma \ref{formule}, hence we have $\varphi_X(x_k(t)n) = \varphi_X(n)$.
\cqfd


\section{Total positivity criteria for $X_K$}
\label{sect8}


Following \cite{GLS4}, lemma 8.3, we put
$N(w) := N \cap w^{-1}N^-w$ and $N'(w) := N \cap w^{-1}N w$, for $w\in W$.
(Obviously, these subgroups do not depend on the choice of a representative of
$w$ in $Norm_G(H)$.)

As shown in \cite{GLS4}, lemma 17.1, we have $N_K = N(w_0w_0^K)$. Similarly, we have :

\begin{lem}\label{ident}
$N \cap B_K^- = N'(w_0w_0^K)$.
\end{lem}

\begin{demo}
Since $w_0^{-1}Nw_0=N^-$, we have
$$ N'(w_0w_0^K) = N \cap ((w_0^K)^{-1}w_0^{-1} N w_0 w_0^K) = N \cap ((w_0^K)^{-1} N^- w_0^K)
=N(w_0^K).$$
We know (see \cite{GLS4}, definition 5.2) that $N(w_0^K)$ is the subgroup of $N$ generated
by the one-parameter subgroups $U_{\alpha}$ for a positive root $\alpha$ such that
$w_0^K(\alpha)$ is a negative root. These are exactly the one-parameter subgroups of $N$
which belong to the Levi subgroup of $B_K^-$, hence $N(w_0^K)=N \cap B_K^-$ and the lemma
follows.
\cqfd

Let us denote $N'(w_0w_0^K)$ by $N_K'$. It is well known that the map
$(n',n) \mapsto n'n$ from $N_K' \times N_K$ to $N$ is a bijection. Hence we have a bijection
between $N_K$ and $N_K' \setminus N$, that is between $N_K$ and $(B_K^- \cap N) \setminus N$ thanks
to lemma \ref{ident}. This bijection coincides with the restriction of $\pi_K$
to $N_K$. Moreover, $\pi_K(N)=\pi_K(N_K)$ and in each fiber
$\pi_K^{-1} \circ \pi_K(x)$ for $x \in N$, there is a unique $n \in N_K$.
All the above maps being regular, we get an isomorphism between
$\ec[N_K]$ and $\ec[N]^{N_K'}:=\ensemble{f \in \ec[N]}{f(n'x)=f(x), \, \forall x \in N, \, \forall n' \in N_K'}$.
Geiss, Leclerc and Schröer have shown (see \cite{GLS3} \S 9 and \cite{GLS4} \S 17) that the
$N_K'$-invariant functions
$\varphi_{M_m}$, where $M_m$ runs over all direct summands of $U_{\textbf{i}}^{(K)}$, form an
initial cluster for a cluster algebra structure on $\ec[N_K] \cong \ec[N]^{N_K'}$.

By definition of $(X_K)_{>0}$, every $y \in (X_K)_{>0}$ belongs to
$\pi_K(N)=\pi_K(N_K)$. Hence it is enough to formulate
our positivity criteria for elements $y$ of $\pi_K(N_K)$.

\begin{thm}[\cite{GLS6}, Conjecture 19.2]\label{resultat1}
Let $y \in X_K$ and suppose that $y=\pi_K(n)$, where $n \in N_K$.
Fix $\textbf{i} \in \RR(w_0,K)$ and let $W=W_1 \oplus \cdots \oplus W_{r_K}$
be a maximal rigid module in $\Sub Q_J$ mutation-equivalent to
$U_{\textbf{i}}^{(K)}$. Then :
$$ y \in (X_K)_{>0} \Longleftrightarrow \left( \forall i=1, \hdots, r_K, \, \varphi_{W_i}(n)>0 \right). $$
\end{thm}

\begin{demo}
Let $y \in \pi_K(N_{>0})=(X_K)_{>0}$. There exists some $x \in N_{>0}$ such that
$y=\pi_K(x)$. Thanks to theorem \ref{Chamber_Ansatz} and lemma \ref{lienV}, we have
$\varphi_{V_k}(x)>0$ for all indecomposable summands $V_k$ of
$V_{m(\textbf{i})}$. But $V_{m(\textbf{i})}$ and $T_{\textbf{i}}^{\dagger}$ are mutation-equivalent
by theorem \ref{algorithme}, hence $\varphi_{M_k}(x)>0$
for all indecomposable summands $M_k$, $k=1,\hdots,r$ of $T_{\textbf{i}}^{\dagger}$.
(Indeed, by definition of the Fomin-Zelevinsky mutations,
each function $\varphi_{M_k}$ is a subtraction-free rational expression in the
functions $\varphi_{V_k}$.)
Thanks to lemma \ref{ident}, there exists some $n \in N_K$
such that $y=\pi_K(n)$. By lemma \ref{invariance}, all the functions $\varphi_U$ for $U$ a summand
of $U_{\textbf{i}}^{(K)}$ belong to $\ec[N]^{N_K'}$, hence,
since $U_{\textbf{i}}^{(K)}$ is a direct summand of $T_{\textbf{i}}^{\dagger}$, $\varphi_{U}(n)=
\varphi_{U}(x)>0$. Thus, if $W=W_1 \oplus \cdots \oplus W_{r_K}$
is mutation equivalent to $U_{\textbf{i}}^{(K)}$ in $\Sub Q_J$, we have $\varphi_{W_j}(n)>0$
for all $1 \leqslant j \leqslant r_K$.

Conversely, when $W=W_1 \oplus \cdots \oplus W_{r_K}$
is mutation equivalent to $U_{\textbf{i}}^{(K)}=U_1 \oplus \cdots \oplus U_{r_K}$ in $\Sub Q_J$, if
$n \in N_K$ is such that $\varphi_{W_i}(n)>0$ for $i=1,\ldots r_K$,
then $\varphi_{U}(n)>0$ for all summands $U$ of $U_{\textbf{i}}^{(K)}$.
Since the application from $N_{>0}$ to $\er^r_{>0}$ which maps $x$ to
$\left( \varphi_{T_1^{\dagger}}(x), \hdots, \varphi_{T_r^{\dagger}}(x) \right)$
is a bijection (where $T_{\textbf{i}}^{\dagger} = T_1^{\dagger} \oplus \cdots
\oplus T_r^{\dagger}$ and the first $r_K$ summands of $T_{\textbf{i}}^{\dagger}$
coincide with the summands of $U_{\textbf{i}}^{(K)}$), there exists an $x \in N_{>0}$ such that
$\varphi_{T_j^{\dagger}}(x)=\varphi_{U_j}(n)>0$ if $j \leqslant r_K$ and $\varphi_{T_j^{\dagger}}(x)=1>0$
if $j>r_K$. Thus $\varphi_T(x)>0$ for all summands $T$ of
${T_{\textbf{i}}^{\dagger}}$, hence $y=\pi_K(x) \in (X_K)_{>0}$.
It only remains to show that $\pi_K(x)=\pi_K(n)$.
Recall that $\left( \varphi_{U_1}, \hdots, \varphi_{U_{r_K}} \right)$ is a seed
for the cluster algebra structure on $\ec[N_K]=\ec[N]^{N_K'}$, thus
every function $f \in \ec[N]^{N_K'}$ is a Laurent polynomial in the variables
$\varphi_{U_j}$ for $1 \leqslant j \leqslant r_K$. But the equations
$\varphi_{T_j^{\dagger}}(x)=\varphi_{U_j}(x)=\varphi_{U_j}(n)$ if $j \leqslant r_K$
imply that for all $f \in \ec[N]^{N_K'}$, we have
$f(n)=f(x)$, that is, $\pi_K(n)=\pi_K(x)$.
\cqfd

Note that we proved the conjecture 19.2 of \cite{GLS6} only in the case where $W=W_1 \oplus \cdots \oplus W_{r_K}$
is a maximal rigid module in $\Sub Q_J$ mutation-equivalent to $U_{\textbf{i}}^{(K)}$ (such a module
is called \emph{reachable}). It remains an open problem to know
whether every maximal rigid module in $\Sub Q_J$ is reachable or not.


\section{Examples}


\subsection{} \label{sect9}
This first example shows that the positivity of all Plücker coordinates is not
a sufficient condition for total positivity. We take $\gg=\mathfrak{sl}_4$ of type $\ea_3$.

Here $I=\{1,2,3\}$. We take $K=\{2\}$, hence $J=\{1,3\}$, and the partial flag variety is :
$$ X_K = \FR_{1,3} := \ensemble{V^1 \subset V^3 \subset \ec^4}{\dim(V^1)=1\,
,\dim(V^3)=3}.$$
We have :
$$N_K = \ensemble{\left( \begin{array}{cccc}
1 & n_{12} & n_{13} & n_{14} \\
0 & 1 & 0 & n_{24} \\
0 & 0 & 1 & n_{34} \\
0 & 0 & 0 & 1 \end{array} \right) }{n_{ij} \in \ec, \, \forall i,j}.$$

Let $\textbf{i} := (2,1,3,2,1,3) \in \RR(w_0,K)$. We have $e(\textbf{i}) = \{1,2,3\}$, $I_K = \{-1,1,-3\}$,
$e_K(\textbf{i}) = \{2,3\}$. The indecomposable (rigid) injectives in $\md \Lambda$ are :
$$\begin{array}{rcccl}
Q_3 &=& M_{-3} &=& \begin{array}{ccc} 1&&\\&2&\\&&3
               \end{array}, \\~\\
Q_2 &=& M_{-2} &=& \begin{array}{ccc} &2&\\1&&3\\&2&
               \end{array}, \\~\\
Q_1 &=& M_{-1} &=& \begin{array}{ccc} &&3\\&2&\\1&&
               \end{array}.
\end{array}$$
Here we use the same convention for representing $\Lambda$-modules as explained in \cite{GLS4} \S 2.4.
Following \cite{GLS3} \S 9.3, we get the three other indecomposable rigid modules in $\Sub Q_J$ (associated to
$\textbf{i}$) :
$$\begin{array}{rcl}
M_1 & = & \begin{array}{ccc} &2& \\ 1&&3 \end{array}, \\~\\
M_2 & = & \begin{array}{ccc} &&3 \end{array}, \\~\\
M_3 & = & \begin{array}{ccc} 1&& \end{array}.
  \end{array}$$
The cluster algebra $\ec[N_K]$ is of type $\ea_1 \times \ea_1$.
We have that $1 \in I_K$, hence $M_1$ is indecomposable injective, and so $\varphi_{M_1}$
is a coefficient in the cluster algebra $\ec[\FR_{1,3}]$,
and thus belongs to all clusters. But the socle
of $M_1$ is not simple, it means that $\varphi_{M_1}$ is not a flag minor. The calculation gives :
$$\begin{array}{rcl}
\varphi_{M_1} & = & D_{s_2 (\varpi_2) , w_0 (\varpi_2)} \\
& = & D_{13 , 34} \\
& = & D_{1,2} D_{123,134} - D_{123,234}.
 \end{array}$$

The other coefficients are $\varphi_{Q_1}=D_{1,4}$ and $\varphi_{Q_3}=D_{123,234}$.
The remaining cluster variables are $\varphi_{M_2}=D_{123,124}$ and $\varphi_{M_3}=D_{1,2}$.
The cluster :
$$\Sigma(U_{\textbf{i}}^{(K)}) = \left\{ \big ( \varphi_{M_2}, \varphi_{M_3},
\varphi_{Q_1}, \varphi_{M_1}, \varphi_{Q_3} \big ) , \Gamma (U_{\textbf{i}}^{(K)}) \right\}$$
is an initial seed for the cluster algebra structure on $\ec[N_K]$. This cluster gives rise to
the following positivity criterium :
$$ (\FR_{1,3})_{>0} = \Ensemble{\pi_K \left( \begin{array}{cccc}
1 & n_{12} & n_{13} & n_{14} \\
0 & 1 & 0 & n_{24} \\
0 & 0 & 1 & n_{34} \\
0 & 0 & 0 & 1 \end{array} \right) }{\begin{array}{rcl}
n_{34} & > & 0 \\
n_{12} & > & 0 \\
n_{14} & > & 0 \\
n_{13}n_{34}-n_{14} &>& 0 \\
n_{14}-n_{13}n_{34}-n_{12}n_{24} &>&0 \end{array}}.$$

Note that the last two inequalities implies $n_{24}<0$.
Note also that if a flag in $\FR_{1,3}$ is totally positive, then its Plücker coordinates
$D_{1,2}$, $D_{1,3}$, $D_{1,4}$, $D_{123,124}$, $D_{123,134}$ and $D_{123,234}$ are
all positive. But the converse is not true. For example :
$$ \left( \begin{array}{cccc}
1&1&1&2 \\ 0&1&0&-1 \\ 0&0&1&1 \\ 0&0&0&1
          \end{array} \right)$$
has all its Plücker coordinates positive but does not belong to $(\FR_{1,3})_{>0}$ since
$n_{13}n_{34}-n_{14} = -1$.


\subsection{}
\label{sect10} \label{sect12}


We take $\gg$ of type $\ed_4$, so $I = \{ 1,2,3,4 \}$.
We label the Dynkin diagram so that the central vertex is $3$.
We will denote by $Q_1$,
$Q_2$, $Q_3$ and $Q_4$ the projective modules in $\md \Lambda$, that is :
$$ \begin{array}{rclcrcl}
Q_1 &=& \begin{array}{ccc}
&1&\\&3&\\2&&4\\&3&\\&1& \end{array}, &~~~~~~~~~~~~&
Q_2 &=& \begin{array}{ccc}
&2&\\&3&\\1&&4\\&3&\\&2& \end{array}, \\~\\
Q_3 &=& \begin{array}{cccc}
&&3&\\1&2&&4\\&&33&\\1&2&&4\\&&3& \end{array}, &~~~~~~~~~~~~&
Q_4 &=& \begin{array}{ccc}
&4&\\&3&\\1&&2\\&3&\\&4&. \end{array}
\end{array}$$

We take $K = \{ 1,2,3 \}$, $J=\{ 4 \}$ and
$\textbf{i} = (1,2,3,1,2,3,4,3,2,1,3,4) \in \RR(w_0,K)$. Here we
have $e(\textbf{i}) = \{ 1,2,3,4,5,6,7,8 \}$, $I_K = \{ -4,4,5,6 \}$ and
$e_K(\textbf{i}) = \{ 7,8 \}$. First we will compute the module $T_{m(\textbf{i})}$.
Note that $m(\textbf{i})=(4,3,1,2,3,4,3,2,1,3,2,1)$. We get :
$$ \begin{array}{rclcrcl}
V_1 &=& 1, &~~~~~~~~~~~~&
V_2 &=& 2, \\~\\
V_3 &=& \begin{array}{ccc}
1&&2\\&3& \end{array}, &~~~~~~~~~~~~&
V_4 &=& \begin{array}{c}
2\\3\\1 \end{array}, \\~\\
V_5 &=& \begin{array}{c} 1\\3\\2 \end{array}, &~~~~~~~~~~~~&
V_6 &=& \begin{array}{ccc}
&3&\\1&&2\\&3& \end{array}, \\~\\
V_7 &=& \begin{array}{ccc}
&3&\\1&&2\\&3&\\&4& \end{array}, &~~~~~~~~~~~~&
V_8 &=& \begin{array}{cccc}
1&2&&\\&&33&\\1&2&&4\\&&3& \end{array}, \end{array}$$
and $V_9=Q_2$, $V_{10}=Q_1$, $V_{11}=Q_3$ and $V_{12}=Q_4$.
Thus the definition of $T_{m(\textbf{i})}$ yields :
$$ T_{m(\textbf{i})} = Q_1 \oplus Q_2 \oplus Q_3 \oplus
\frac{Q_1}{V_1} \oplus
\frac{Q_2}{V_2} \oplus
\frac{Q_3}{V_3} \oplus
Q_4 \oplus \frac{Q_3}{V_6} \oplus \frac{Q_2}{V_5}
\oplus \frac{Q_1}{V_4} \oplus \frac{Q_3}{V_8} \oplus \frac{Q_4}{V_7}.$$
Then, we compute $T_{\textbf{i}}^{\dagger}$ :
$$ T_{\textbf{i}}^{\dagger} = \bigoplus_{i \in I} Q_i \oplus \bigoplus_{m \in
e(\textbf{i})} N_m,$$
where
$$ \begin{array}{rclcrcl}
N_1 &=& \begin{array}{ccc}
&1&\\&3&\\2&&4\\&3& \end{array}, &~~~~~~~~~~~~&
N_2 &=& \begin{array}{ccc}
&2&\\&3&\\1&&4\\&3& \end{array}, \\~\\
N_3 &=& \begin{array}{cccc}
&&3&\\1&2&&4\\&&33&\\&&&4 \end{array}, &~~~~~~~~~~~~&
N_4 &=& \begin{array}{c} 1\\3\\4 \end{array}, \\~\\
N_5 &=& \begin{array}{c} 2\\3\\4 \end{array}, &~~~~~~~~~~~~&
N_6 &=& \begin{array}{cccc}
&&3&\\1&2&&4\\&&3&\\&&&4 \end{array}, \\~\\
N_7 &=& 4, &~~~~~~~~~~~~&
N_8 &=& \begin{array}{cc} 3\\4 \end{array}. \end{array}$$
We can see that :
$$ T_{m(\textbf{i})} = Q_1 \oplus Q_2 \oplus Q_3 \oplus N_1 \oplus N_2 \oplus N_3 \oplus Q_4 \oplus N_6 \oplus N_5
\oplus N_4 \oplus N_8 \oplus N_7,$$
that is, $T_{m(\textbf{i})} = T_{\textbf{i}}^{\dagger}$,
which illustrates Lemma \ref{lien}.

Here, $ U_{\textbf{i}}^{(K)} = N_4 \oplus N_5 \oplus N_6 \oplus Q_4 \bigoplus
N_7 \oplus N_8 $.
The cluster algebra $\ec[N_K]$ is of finite type $\ea_1 \times \ea_1$.
Note that the socle of $N_6$ is isomorphic to $S_4 \oplus S_4$, thus is not
simple, but $N_6$ is projective in $\Sub Q_4$, hence $\varphi_{N_6}$ is a
coefficient of the cluster algebra $\ec[N_K]$ which is not a flag minor.
The module $U_{\textbf{i}}^{(K)}$ gives the following
positivity criterium (thanks to Theorem \ref{resultat1}) :
$$(X_K)_{>0}= \Ensemble{\pi_K \left( n \right) }{\begin{array}{ccc}
\varphi_{N_4}(n)>0, &\,& \varphi_{N_5}(n)>0,\\
\varphi_{N_6}(n)>0,&&
\varphi_{N_7}(n)>0,\\
\varphi_{N_8}(n)>0,&&
\varphi_{Q_4}(n)>0 \end{array}}.$$
Recall that one can realize every irreducible highest weight $\gg$-module $L(\lambda)$
as a subspace of $\ec[N]$ (see \cite{GLS6}, \S 8). The dual of Lusztig's
semicanonical basis for $L(\varpi_4)$ (seen as a subspace of $\ec[N]$) consists on
the functions $\varphi_{N_7}$, $\varphi_{N_8}$, $\varphi_{N_5}$, $\varphi_{N_4}$,
$\varphi_{Y}$, $\varphi_{V_7}$, $\varphi_{Q_4}$ and the constant function $1$
(these eight functions are the Plücker coordinates of the flag variety), where
$Y$ is the following $\Sub Q_4$-module :
$$ Y = \begin{array}{ccc}
1&&2\\&3&\\&4& \end{array}.$$
The cluster algebra $\ec[N_K]$ comes
with only two mutation relations : $\varphi_{N_7} \varphi_{V_7} = \varphi_{N_6} + \varphi_{Q_4}$
and $\varphi_{N_8} \varphi_{Y} = \varphi_{N_6} + \varphi_{N_4} \varphi_{N_5}$.
We see that if $x \in (X_K)_{>0}$ and $n \in N_K$ is such that $\pi_K(n)=x$, then
$\varphi_{Y}(n)>0$ and $\varphi_{V_7}(n)>0$ (thanks to the preceding mutation relations), hence
all the Plücker coordinates of $n$ are positive. But the converse is not true,
since the positivity of these Plücker coordinates do not imply the positivity of $\varphi_{N_6}$
which is needed in all positivity criteria (since $\varphi_{N_6}$ is a coefficient
of the cluster algebra). Since $\varphi_{N_6}$ is not in $L(\varpi_4)$ but is in
$L(2 \varpi_4)$, the hypothesis of \cite{Lu1} theorem 3.4 can not be weakened,
in contradiction with what is stated in the note of \cite{Lu1} \S 3.12.

\bigskip
\small
\noindent
\begin{tabular}{ll}
Nicolas {\sc Chevalier} : &
LMNO, CNRS UMR 6139,
Universit\'e de Caen,\\
& F-14032 Caen cedex, France\\
&email : {\tt nicolas.chevalier01@unicaen.fr}
\end{tabular}

\end{document}